\documentclass{roadef_en}

\usepackage[backend=biber, maxbibnames=99]{biblatex} 
\usepackage{hyperref}
\usepackage{graphicx}
\usepackage{fancyhdr}

\fancypagestyle{titlelogos}{%
  \fancyhf{}
  \fancyhead[L]{%
    \includegraphics[height=1.2cm]{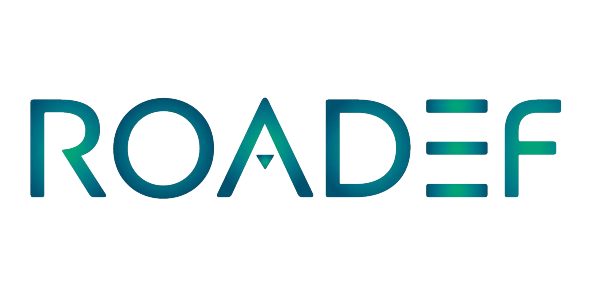}\hspace{0.5cm}%
    \includegraphics[height=1.2cm]{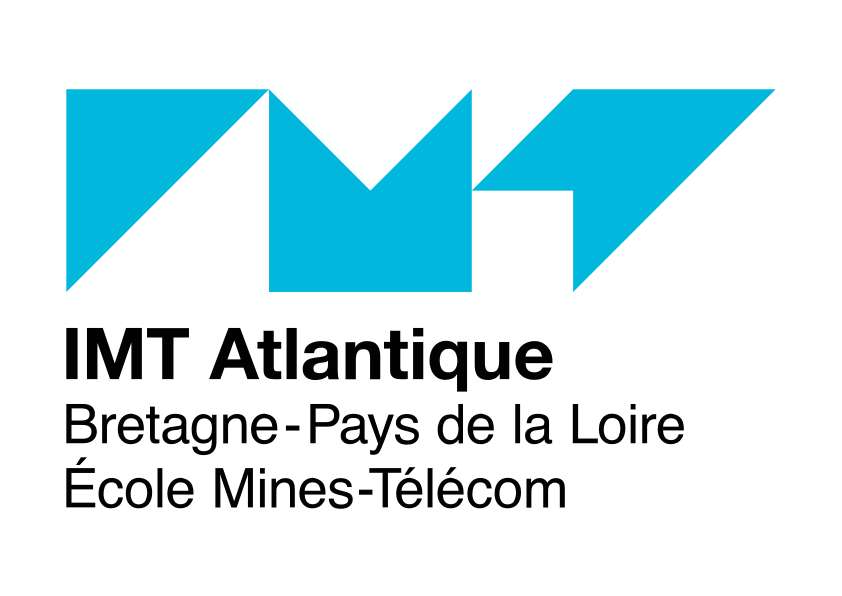}\hspace{0.5cm}%
    \includegraphics[height=1.2cm]{logo_insa_lyon.jpg}\hspace{0.5cm}%
    \includegraphics[height=1.2cm]{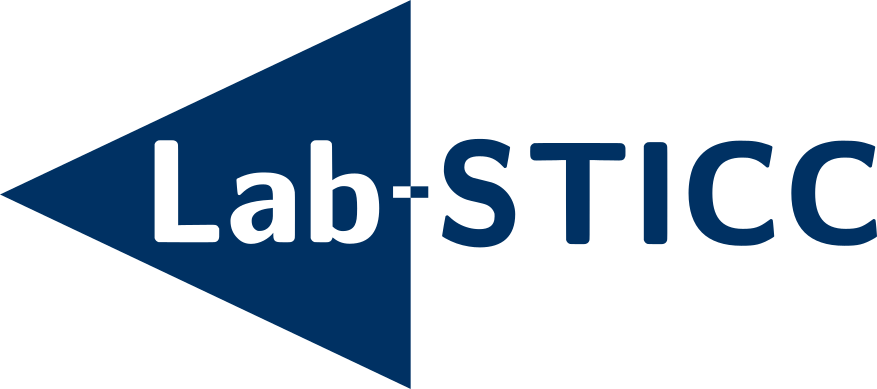}\hspace{0.5cm}%
    \includegraphics[height=1.2cm]{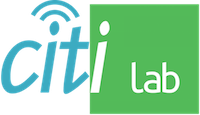}%
  }%
}

\begin{document}

\title{Randomized Constructive Heuristics for the VRPTW:\\ A Focus on Regret-k}

\def\shorttitle{Randomized Constructive Heuristics for the VRPTW}

\author{
    Florian Rascoussier\inst{1, 2}, 
    Romain Billot\inst{1},
    Lina Fahed\inst{1},
    Christine Solnon\inst{2}
}

\institute{
IMT Atlantique, Lab-STICC, UMR CNRS 6285, 29238 Brest, France  \\
\email{florian.rascoussier@imt-atlantique.fr}
\and
INSA Lyon, Inria, CITI, UR3720, 69621 Villeurbanne, France \\
\email{florian.rascoussier@insa-lyon.fr}
}

\maketitle
\thispagestyle{titlelogos}  


\section{Introduction}

The Vehicle Routing Problem with Time Windows (VRPTW) is a classic combinatorial optimization problem in which a fleet of capacitated vehicles must serve a set of customers within strict time windows while minimizing the total travel time. Constructive heuristics are fundamental for generating solutions, as seen in Ant Colony Optimization (ACO), and often serve as starting points for metaheuristics such as Iterated Local Search (ILS) or Genetic Algorithms (GA). This work focuses on three distinct constructive heuristics: the Nearest Neighbor (NN), Best Insertion (BI), and Regret-k heuristics. NN builds routes sequentially, extending the current route by visiting the closest feasible unvisited customer. BI constructs routes in parallel\footnote{Note that the term \textit{parallel} refers to the construction procedure and not necessarily to parallel computing, as coined in the heuristic literature \cite{Solomon1987, potvinParallelRouteBuilding1993}}, inserting the unvisited customer that minimizes the cost increase into the best position across all routes \cite{Solomon1987}. The Regret-k heuristic 
extends BI by incorporating a look-ahead mechanism \cite{potvinParallelRouteBuilding1993}: it calculates a regret value for each customer, defined as the cost difference between its best insertion and its $k$-th best insertion \cite{foisyImplementingInsertionHeuristic1993}. By prioritizing customers with high regret---those that would become significantly more expensive to insert later---it is generally considered as one of the most effective greedy heuristic. 


Although these heuristics are well-established as deterministic initializers, their potential within a randomized, multi-start framework remains largely under-explored. We investigate three primary research questions: (1) How can we effectively randomize these heuristics to implement a pure Multi-start Greedy Randomized (GR) approach? (2) Specifically, how can the Regret-k heuristic be randomized? 
(3) What is the impact of 
combining these greedy randomized constructions with LS and Ant Colony Optimization (ACO).

\section{Methodology}


Similarly to ACO-based metaheuristics, we handle NN randomization using probabilities pondered by the heuristic values of each evaluated feasible extension; a principle that can be directly transposed to BI. To address the randomization of the Regret-k heuristic, we distinguish between the ``who'' and ``where'' decisions. The ``who'' decision concerns which customer to insert next, while the ``where'' decision concerns its insertion position. In its standard deterministic form, Regret-k selects the customer with the maximum regret value. We propose to randomize this selection process by choosing the next customer using a probabilistic selection proportional to the regret. This approach preserves the heuristic's foresight while introducing the necessary diversity for a multi-start framework. We also explore the impact of randomizing the ``where'' decision, as well as the combination of both randomizations. We then study the impact of LS for post-construction improvement before exploring the potential for pheromonal online learning in ACO settings. 


\section{Experimentation}

We evaluate the performance of the proposed randomized heuristics on the classic Solomon \cite{Solomon1987} and Gehring \& Homberger \cite{HombergerGehring1999} benchmark instances, ranging from $100$ to $1000$ customers. To ensure fair and reproducible comparisons, we strictly adhere to the DIMACS 2021 competition conventions \cite{DIMACS2021}, utilizing integer precision and truncated Euclidean distances. The results are compared against the Best Known Solutions (BKS) and the state-of-the-art Hybrid Genetic Search (HGS) solver \cite{PyVRP2024, vidalHybridGeneticSearch2022}, reusing their LS implementation. All critical heuristic components are implemented in modern C++ and bound to Python using pybind11, following the high-performance architecture of the PyVRP project.


\section{Conclusion}

This study provides a comprehensive analysis of randomized constructive heuristics for the VRPTW, with a particular emphasis on the novel randomization of the Regret-k heuristic. To the best of our knowledge, although Regret-k is widely used for (TD)-VRPTW initialization, its randomization has never been examined. 
By separating customer-selection and insertion decisions, we show how diversity can be injected into this powerful heuristic. Our extensive benchmarking on standard datasets highlights the trade-off between solution quality and computational cost, offering insights for the design of efficient metaheuristics.

This work constitutes a preliminary step towards more advanced hybrid methods, and we plan to integrate these findings into ACO, ILS and Branch-and-Price algorithms in future research on the Time-Dependent VRPTW (TDVRPTW) as part of the \href{https://anr.fr/Project-ANR-22-CE22-0016}{MAMUT} project.

\printbibliography
\end{document}